\documentclass[11pt]{article}
\usepackage{amscd}
\usepackage{amsmath}
\usepackage{amsfonts}

\newcommand{\ep}{\varepsilon}

\newcommand{\ddbar}{\frac{\sqrt{-1}}{2\pi} \partial \overline{\partial}}

\newcommand{\ov}[1]{\overline{#1}}

\newcommand{\oke}{\omega_{\textrm{KE}}}
\newcommand{\hke}{h_{\textrm{KE}}}
\newcommand{\Xcan}{X_{\textrm{can}}}
\newcommand{\Xmin}{X_{\textrm{min}}}
\newcommand{\hmin}{h_{\textrm{min}}}
\newcommand{\omin}{\omega_{\textrm{min}}}
\newcommand{\HZ}{\textrm{H}^0}
\newcommand{\SE}{S_{-1}}
\newcommand{\SM}{S_{-2}}
\newcommand{\hC}{h_{\mathcal{C}}}

\begin{document}
\newcounter{remark}
\newcounter{theor}
\setcounter{remark}{0}
\setcounter{theor}{1}
\newtheorem{claim}{Claim}
\newtheorem{theorem}{Theorem}[section]
\newtheorem{proposition}{Proposition}[section]
\newtheorem{lemma}{Lemma}[section]
\newtheorem{defn}{Definition}[theor]
\newtheorem{corollary}{Corollary}[section]
\newenvironment{proof}[1][Proof]{\begin{trivlist}
\item[\hskip \labelsep {\bfseries #1}]}{\end{trivlist}}
\newenvironment{remark}[1][Remark]{\addtocounter{remark}{1} \begin{trivlist}
\item[\hskip
\labelsep {\bfseries #1  \thesection.\theremark}]}{\end{trivlist}}
\newenvironment{proof1}{\noindent\textbf{Proof} }{\hfill$\square$\medskip}

\centerline{\bf CONSTRUCTIONS OF K\"AHLER-EINSTEIN METRICS}
\centerline{\bf WITH NEGATIVE SCALAR CURVATURE}

\bigskip
\begin{center}
\begin{tabular}{ccc}
{\bf Jian Song\footnotemark[1]}  & & {\bf Ben Weinkove\footnotemark[2]} 
 \\
Johns Hopkins University & & Harvard University \\
Department of Mathematics & & Department of Mathematics \\
Baltimore, MD 21218 & & Cambridge, MA 02138 \\
\end{tabular}
\end{center}
\footnotetext[1]{The first-named author is on leave for the semester and is visiting MSRI, Berkeley, CA as a postdoctoral fellow.   He is supported in part by National Science Foundation grant DMS 0604805.}
\footnotetext[2]{Part of this research was carried out while the second-named author was a short-term visitor at MSRI in January 2007.  He is supported in part by National Science Foundation grant DMS 0504285.}

\bigskip
\bigskip
\noindent
{\bf Abstract.} \ We show that on K\"ahler manifolds with negative first Chern class, the
sequence of algebraic metrics introduced by H. Tsuji converges uniformly
to the K\"ahler-Einstein metric.  For algebraic surfaces of general type and
orbifolds with isolated singularities, we prove a convergence result for a
modified version of Tsuji's iterative construction.

\setcounter{equation}{0}
\setcounter{lemma}{0}
\addtocounter{section}{1}
\setlength{\arraycolsep}{2pt}

\bigskip
\bigskip
\noindent
{\bf 1.  Introduction}
\bigskip

A K\"ahler manifold with negative first Chern class admits a unique K\"ahler-Einstein metric.  This was shown by Yau in his seminal paper on the Calabi conjecture \cite{Y1}, and also independently by Aubin \cite{A}.  Yau later posed the question of whether, in general, K\"ahler-Einstein metrics could be obtained as a limit of algebraic metrics induced from embeddings into projective space \cite{Y2}.   Donaldson \cite{D1} showed that if a polarized variety $(X, L)$ with discrete automorphism group admits a constant scalar curvature K\"ahler metric, then there is indeed a sequence of algebraic `balanced metrics' \cite{Zh} converging to it.  Donaldson's proof makes use of the Tian-Yau-Zelditch expansion of the Bergman kernel \cite{Ti}, \cite{Ze} (see also \cite{C}) and Lu's \cite{L} computation of the second coefficient.  Very recently, Donaldson \cite{D2} has described how numerical approximations to these balanced metrics could be used to compute, to a good accuracy, explicit K\"ahler-Einstein metrics on certain varieties.  

Tsuji \cite{Ts} has considered a different way of producing K\"ahler-Einstein metrics by algebraic approximations.  He introduced a new iterative procedure on varieties of general type with the aim of describing (possibly singular) K\"ahler-Einstein metrics.    In the case when the first Chern class is negative, Tsuji proved that his iteration converges, in a certain weak sense, to the K\"ahler-Einstein metric.  In this paper, we give a uniform convergence result and describe how his procedure may be modified to obtain results in the case of algebraic surfaces of general type, and on orbifolds with isolated singuarities.  

We now describe Tsuji's iteration.  Let $X$ be a compact K\"ahler manifold of complex dimension $n$ with  ample canonical bundle $K_X$.  Let $\oke = \frac{\sqrt{-1}}{2\pi} (g_{\textrm{KE}})_{i \ov{j}} dz^i \wedge dz^{\ov{j}}$ be the K\"ahler-Einstein metric satisfying 
$$ \textrm{Ric}(\oke) = - \frac{\sqrt{-1}}{2\pi} \partial \overline{\partial} \log \oke^n = - \oke \in c_1(X).$$
Fix a Hermitian metric $\hke$ on $K_X$ by setting $\hke = (\det g_{\textrm{KE}})^{-1}$.  

Let $m_0\ge 1$ be an integer such that $K_X^{m_0}$ is base point free and let $h_{m_0}$ be any Hermitian metric on $K_X^{m_0}$.  Define a sequence of Hermitian metrics $h_m$ on $K_X^m$ for $m > m_0$ inductively as follows.  Suppose that $h_m$ is a given Hermitian metric on $K_X^m$.  To define $h_{m+1}$,  first define an inner product $\langle \ , \ \rangle_{T_{m+1}}$ on the space of sections of $K_X^{m+1}$ by
\begin{equation} \label{eqnT}
\langle s, t \rangle_{T_{m+1}} = \int_X h_m \otimes s \otimes \ov{t},\quad \textrm{for } s, t \in \HZ(X, K_X^{m+1}),
\end{equation}
where $h_m \otimes s \otimes \ov{t}$ is regarded as a volume form on $X$.  Now let $( \sigma_{m+1}^{(0)}, \ldots, \sigma_{m+1}^{(N_{m+1})})$, for $N_{m+1} +1 = \textrm{dim} \, \HZ(X, K_X^{m+1})$, be an orthonormal basis of $\HZ(X, K_X^{m+1})$ with respect to this inner product.
Then define the Hermitian metric $h_{m+1}$ on $K_X^{m+1}$ by
$$h_{m+1} =  \frac{(m+n+1)!}{(m+1)!} \left( \sum_{i=0}^{N_{m+1}}  \sigma_{m+1}^{(i)} \otimes  \ov{\sigma_{m+1}^{(i)}} \right)^{-1}.$$
Observe that this metric is independent of the choice of orthonormal basis.  

We have the following theorem on the convergence of the metrics $h_m$, strengthening the result given in \cite{Ts}.  

\bigskip
\noindent
{\bf Theorem 1}  \emph{Let $X$ be a compact K\"ahler manifold with ample canonical bundle.  Let $h_{\emph{KE}}$, $\omega_{\textrm{KE}}$ and the sequence of Hermitian metrics $h_m$ be as above.  There exists a constant $C$ depending only on $X$ and $h_{m_0}$ such that for $m\ge m_0$, \begin{equation} \label{mainestimate}
e^{-\frac{ C \log m}{m}} h_{\emph{KE}} \le h_m^{1/m} \le e^{\frac{C \log m}{m}} h_{\emph{KE}}.
\end{equation}
Hence $h_m^{1/m}$ converges uniformly on $X$ to $h_{\emph{KE}}$ as $m \rightarrow \infty$.}
\bigskip

We now describe a modification of Tsuji's iteration.  Let $\beta$ be a continuous function on a variety $X$ with $0 \le \beta \le 1$.  Here we allow $X$ to have mild singularities.  We also remove the assumption that $K_X$ be ample.  We only require that there exists an $m_0 \ge 1$ such that $K_X^{m_0}$ is base point free.  Let $h_{m_0}$ be a Hermitian metric on $K_X^{m_0}$.  Given $0<\ep \le 1$, we inductively define a sequence of Hermitian metrics $h_{m, \ep} = h_{m, \ep}(\beta, h_{m_0})$ on $K_X$ as follows.
Assuming that $h_{m, \ep}$ is given, define an inner product $\langle \ , \ \rangle_{T_{m+1, \ep}}$ on the space of sections of $K_X^{m+1}$ by
$$ \langle s, t \rangle_{T_{m+1, \ep}} = \int_X \beta^{\ep} h_{m,\ep} \otimes s \otimes \ov{t}, \quad \textrm{for} \ s,t \in \textrm{H}^0(X, K_X^{m+1}).$$  
Then define the Hermitian metric $h_{m+1, \ep}$ on $K_X^{m+1}$ by
$$h_{m+1, \ep} =  \frac{(m+n+1)!}{(m+1)!} \left( \sum_{i=0}^{N_{m+1}}  \sigma_{m+1, \ep}^{(i)} \otimes  \ov{\sigma_{m+1, \ep}^{(i)}} \right)^{-1},$$
where $( \sigma_{m+1, \ep}^{(0)}, \ldots, \sigma_{m+1, \ep}^{(N_{m+1})})$ is an orthonormal basis of $\HZ(X, K_X^{m+1})$ with respect to the inner product $\langle \ , \ \rangle_{T_{m+1, \ep}}$.  We call $h_{m, \ep}$ the \emph{modified Tsuji iteration}.  It depends on $\ep$ and $\beta$.

First, we consider the case when $X$ is an algebraic surface of general type.  Let $\mathcal{E} = \sum_i E_i$ be the sum of the nonsingular rational curves $E_i$ of self intersection $-1$ (or $(-1)$-curves, for short)  on $X$.
Let $\tau: X \rightarrow \Xmin$ be a holomorphic map blowing down these curves on $X$, so that $\Xmin$ is a minimal surface of general type.   Now if $h$ is any Hermitian metric on $K_{\Xmin}$, then $\Omega = h^{-1}(z^1, \ldots, z^n) (\sqrt{-1}/2\pi)^n dz^1 \wedge d\overline{z^1} \wedge \cdots \wedge dz^n \wedge d\overline{z^n}$ can be regarded as a volume form on $\Xmin$.  Using
 coordinates  $w^i$ on $X$, there is a holomorphic section $\SE$ of the line bundle $[\mathcal{E}]$ associated to $\mathcal{E}$ and vanishing of order one on $\mathcal{E}$ satisfying
\begin{eqnarray*}
\left( (\tau^*h)^{-1}\otimes |\SE|^2 \right)(w^1, \ldots, w^n)  \left( \frac{\sqrt{-1}}{2\pi} \right)^n dw^1 \wedge d\overline{w^1} \wedge \cdots \wedge dw^n \wedge d\overline{w^n} &=&  \tau^*\Omega,
\end{eqnarray*}
for any such $h$.

$\Xmin$ has no $(-1)$-curves, but may have $(-2)$-curves.  Let $f : \Xmin \rightarrow \Xcan$ be its canonical map.  $\Xcan$ is an surface with ample canonical bundle $K_{\Xcan}$ and, at worst, isolated orbifold singularities.  By the orbifold version of the results of \cite{Y1}, \cite{A} (see \cite{K}, for example), there exists an 
 orbifold K\"ahler-Einstein metric   $\omega_{\textrm{KE}}$ on $\Xcan$, with corresponding Hermitian metric $h_{\textrm{KE}}$ on $K_{\Xcan}$. Define a K\"ahler metric on $\Xmin$ by $\omin = f^* \omega_{\textrm{KE}}$ and a Hermitian metric on $K_{\Xmin}$ by $\hmin = f^* h_{\textrm{KE}}$.  Note that $\omin$ and $\hmin$ are not smooth in general, although $\hmin$ is continuous on $\Xmin$ (see \cite{TZ}).  Let $\mathcal{C} = \sum_i C_i$ be the sum of the $(-2)$-curves on $\Xmin$ and let $\SM$ be a holomorphic section of the line bundle $[\mathcal{C}]$, vanishing of order one on $\mathcal{C}$.   Fix a smooth Hermitian metric $h_{\mathcal{C}}$ on $[\mathcal{C}]$, and assume that $\sup_{\Xmin} |\SM|^2_{\hC} = 1$.  Let $\beta$ be the smooth function on $X$ defined by $\beta = \tau^* |\SM|^2_{\hC}$, and let $h_{\infty} = \tau^* \hmin$.

For this $\beta$ and some initial Hermitian metric $h_{m_0}$ on $K_X^{m_0}$,  let $h_{m, \ep} = h_{m, \ep}(\beta, h_{m_0})$ be the sequence of Hermitian metrics in the modified Tsuji iteration as described above.  Then we have the following result.

\bigskip \pagebreak[3]
\noindent
{\bf Theorem 2} \emph{Let $X$ be an algebraic surface of general type.  With $h_{m, \ep} = h_{m, \ep}(\beta, h_{m_0})$ and $\SE$ as described above, for every sequence $\ep_j \rightarrow 0$,
$$\limsup_{m \rightarrow \infty} h_{m, \ep_j}^{1/m} \rightarrow h_{\infty} \otimes |\SE|^{-2}, \qquad \textrm{as } j \rightarrow \infty,$$
almost everywhere on $X$.}
\bigskip

Note that since $K_X = \tau^* K_{\Xmin} + [\mathcal{E}]$, we can regard $h_{\infty} \otimes |\SE|^{-2}$  as a singular Hermitian metric on $K_X$.  

We now consider the case when $X$ is an orbifold with isolated singularities with ample canonical bundle.

\bigskip
\noindent
{\bf Theorem 3} \emph{Let $(X, \omega_{\emph{KE}})$ be a K\"ahler-Einstein orbifold with $K_X$ ample and only isolated singularities at points $p_1, \ldots, p_k$.  Let $\beta$ be a continuous function on $X$ satisfying $0 \le \beta \le 1$ and $\beta(x) = 0$ if and only if $x=p_i$ for some $i$.  Then, with $h_{m, \ep} = h_{m, \ep}(\beta, h_{m_0})$ as above, for every sequence $\ep_j \rightarrow 0$,
$$\limsup_{m \rightarrow \infty} h_{m, \ep_j}^{1/m} \rightarrow \hke, \quad \textrm{as } j \rightarrow \infty,$$
almost everywhere on $X$.}
\bigskip

\noindent
We end the introduction with a couple of  remarks.
\begin{enumerate}
\item  Tsuji's iteration has some similarities to Donaldson's $T_K$-iteration (see \cite{D2}, section 2.2.2) which in the case of $K_X$ ample should yield a `canonically balanced metric' using sections of a fixed power $k$ of the canonical line bundle.  As $k \rightarrow \infty$, the limit of these metrics is expected to be the K\"ahler-Einstein metric.  On the other hand, Tsuji's method is a single iterative process.  
\item  Donaldson has suggested that the dynamical systems introduced in \cite{D2} are likely to be discrete approximations to the K\"ahler-Ricci and Calabi flows.  It would be interesting to know whether Tsuji's iteration could be viewed in a similar light.
\end{enumerate}

The outline of the paper is as follows.  Our main technique is the peak section method of \cite{Ti}.  This is described in section 2, and extended to orbifolds with isolated singularities (for related results on the Szeg\"o kernel for orbifolds, see \cite{S}, \cite{DLM}, \cite{P}).  In sections 3 and 4 we prove the main theorems.

\setcounter{equation}{0}
\setcounter{lemma}{0}
\addtocounter{section}{1}
\bigskip
\bigskip
\noindent
{\bf 2. Peak sections}
\bigskip

Now suppose that $(X, \omega)$ be a compact K\"ahler manifold of complex dimension $n$ with $\omega \in c_1(L)$ for an ample line bundle $L$.  
Fix a Hermitian metric $h$ on $L$ satisfying $- \ddbar \log h = \omega$.   Write $\langle  \cdot  , \cdot \rangle_{L^2(h^m)}$ and $\| \cdot \|_{L^2(h^m)}$ for the $L^2$ inner product and norm on $\HZ(X, L^m)$ with respect to the Hermitian metric $h^m$ and volume form $\frac{1}{n!} \omega^n$.
We use the following lemma from \cite{Ti}.

\pagebreak[3]
\begin{lemma} \label{lemmapeaksections}
There exists $m_1>0$ depending only on $X$, $L$ and $h$ such that for every $x_0 \in X$ and $m \ge m_1$ there is a global holomorphic section $s_{m, x_0}$ of $L^m$ satisfying the following.
\begin{enumerate}
\item[(i)] $\displaystyle{|s_{m, x_0}|^2_{h^m} (x_0)  =  1}$ and 
\begin{eqnarray*}
%|s_m|^2_{h^m} (x_0) & = & 1, \\
 \| s_{m, x_0} \|^2_{L^2(h^m)} & = & \frac{m!}{(m+n)!} (1 + \emph{O}(m^{-1})),
 \end{eqnarray*}
where $f(m) = \emph{O}(m^{-1})$ means that $|f(m)| m \le A$ for a constant $A$ depending only on $X$, $L$ and $h$.
\item[(ii)] For $t$ any holomorphic section of $L^m$ which vanishes at $x_0$,
$$\left| \langle s_{m, x_0}, t \rangle_{L^2(h^m)} \right| \le \frac{B}{m} \| s_{m, x_0} \|_{L^2(h^m)} \| t \|_{L^2(h^m)},$$
with a constant $B$ depending only on $X$, $L$ and $h$.
\end{enumerate}
\end{lemma}
\begin{proof1}  We outline here the proof of part (i), since we will explicitly refer to this method in the orbifold case.   For (ii) we refer the reader to \cite{Ti}. 
Pick a normal coordinate chart $(U, ( z^1, z^2, \ldots, z^n) )$ for $g$ centered at the point $x_0$.
Let $\eta: [0,\infty) \rightarrow [0,1]$ be a cut-off function satisfying $\eta(t)  =  1$ for $t \le \frac{1}{2}$, $\eta(t)  =  0$ for $t \ge 1$
%$$\begin{array}{rcll} \eta(t) & = & 1,& \qquad t \le \frac{1}{2} \\
%\eta(t) & = & 0, & \qquad t \ge 1, \end{array}$$
and $-4 \le \eta'(t) \le 0$, $| \eta''(t)| \le 8$.
Define a weight function $\psi$, which for $m$ sufficiently large is
supported in $U$, by setting
\begin{equation} \label{eqnpsi}
\psi(z) =  (n+2) \eta\left( a_m |z|^2 \right) \log \left( a_m |z|^2 \right),
\end{equation}
for $z \in U$, and $\psi \equiv 0$ outside $U$, where
$a_m = m/(\log m)^2.$
A short calculation shows that
%\begin{eqnarray*}
%\ddbar \psi (z)  
% & = & (n+2) \left\{ (\log ( a_m |z|^2 ))  \eta''(a_m |z|^2 ) a_m^2 \sqrt{-1} \partial |z|^2 \wedge \ov{\partial} |z|^2 \right.  \\
%&& \mbox{} +  (\log ( a_m |z|^2 )) \eta'(a_m |z|^2) a_m \ddbar |z|^2  \\
%&& \mbox{} + 2 \textrm{Re} \left( \eta'(a_m |z|^2) a_m \partial \log|z|^2 \wedge \ov{\partial} |z|^2 \right) \\
%&& \mbox{} \left. + \eta(a_m |z|^2) \ddbar \log |z|^2 \right\}.
%\end{eqnarray*}
%We now claim that
\begin{equation} \label{eqnddbarpsi}
\ddbar \psi(z) \ge - C (n+2) a_m  \omega, \quad \textrm{for} \  |z| >0. 
\end{equation}
%To see (\ref{eqnddbarpsi}), observe first that if $\eta''(a_m |z|^2) \neq 0$ then $|z|^2$ lies in the interval $[1/(2a_m), 1/a_m]$.  Then
%$$\left| \sqrt{-1} \partial |z|^2 \wedge \ov{\partial} |z|^2 \right|_g \le \frac{C}{a_m},$$
%and 
%$$ - \log 2 \le \log (a_m |z|^2) \le 0.$$
%Hence
%\begin{equation} \label{eqnineq1}
%\left| (\log (a_m |z|^2)) \eta''(a_m |z|^2) a_m^2 \sqrt{-1} \partial |z|^2 \wedge \ov{\partial} |z|^2 \right|_g \le C a_m.
%\end{equation}
%Moreover, since $\eta' \le 0$, we obtain
%\begin{equation} \label{eqnineq2}
%( \log (a_m |z|^2) ) \eta'(a_m |z|^2) \ddbar |z|^2 \ge 0.
%\end{equation}
%Also, we have
%\begin{equation} \label{eqnineq3}
%\left| \eta'(a_m |z|^2) a_m \partial \log |z|^2 \wedge \ov{\partial} |z|^2 \right|_g \le C a_m.
%\end{equation}
%Combining (\ref{eqnineq1}), (\ref{eqnineq2}) and (\ref{eqnineq3}), and using the fact that $\ddbar \log |z|^2=0$ for $z\neq 0$ we obtain (\ref{eqnddbarpsi}).
Now let $\psi_i$ be a decreasing sequence of smooth functions on $X$ such that
\begin{equation} \label{eqnpsii1}
\psi = \lim_{i \rightarrow \infty} \psi_i \quad \textrm{and} \quad  \ddbar \psi_i\ge - C (n+2) a_m  \omega,
\end{equation}
where the constant $C$ may be different from the one in (\ref{eqnddbarpsi}).  It follows that for sufficiently large $m$, 
\begin{equation} \label{eqnpsii2}
\ddbar \psi_i - \ddbar \log h^m + \textrm{Ric}(\omega) \ge \frac{m}{C} \omega.
\end{equation}
Choose a local holomorphic section $v$ of $L$ so that $|v|^2_{h} (x_0) = 1$ and $(\partial | v|^2_{h}) (x_0) = 0$.  Let $w$ be the smooth local section of $L^m$ defined by
$$w = \ov{\partial} \left( \eta \left( \frac{a_m |z|^2}{4} \right) \right) v^m.$$
We apply the $L^2$ estimate of H\"ormander \cite{H} (cf. Proposition 1.1 of \cite{Ti})  to the (1,0) form $w$:  
%\begin{proposition} \label{prophormander} Let $\psi$ be a function on $X$ which can be approximated %by a decreasing sequence of smooth functions $\{ \psi_i \}$.  If
%$$\ddbar \psi_i - \ddbar \log h + \emph{Ric}(\omega) \ge A \omega$$
%for some positive constant $A$ then for any smooth $L$-valued (0,1)-form $w$ on $X$ with $\ov{\partial} w=0$ and
%$\int_X |w|_{h,g}^2 e^{-\psi} \frac{\omega^n}{n!} < \infty,$
%there exists a global smooth section $u$ of $L$ with $\ov{\partial} u =w$ and
%$$\int_X |u|^2_h e^{-\psi} \frac{\omega^n}{n!} \le \frac{1}{A} \int_X | w |^2_{h,g} e^{-\psi} \frac{\omega^n}{n!}.$$
%\end{proposition}
there exists a smooth global section $u$ of $L^m$ satisfying $\ov{\partial} u =w$ and
$$ \int_X |u|^2_{h^m} e^{-\psi} \frac{\omega^n}{n!} \le \frac{C}{m} \int_X |w|_{h^m}^2 e^{-\psi} \frac{\omega^n}{n!}.$$
Here we are using the fact that the weight function $\psi$ can be approximated by $\psi_i$ satisfying  (\ref{eqnpsii2}).  Observe that $w$ vanishes identically outside the region $2/a_m \le |z|^2 \le 4/a_m$, and that in $U$,
$| w|^2_{h^m} e^{-\psi} \le C a_m |v^m|^2_{h^m}.$
It follows that
$$ \int_X |u|^2_{h^m} e^{-\psi} \frac{\omega^n}{n!} \le \frac{C}{(\log m)^2} \int_{2/a_m \le |z|^2 \le 4/a_m} |v^m|^2_{h^m} \frac{\omega^n}{n!}.$$
From our choice of coordinates,
$|v|^2_{h}(z) = 1 - |z|^2 + \textrm{O}(|z|^3),$
and so locally
\begin{eqnarray*}
|v^m|^2_{h^m} \frac{\omega^n}{n!} 
& \le &  \left(1- \frac{1}{2}|z|^2 \right)^{m} \left(\frac{\sqrt{-1}}{2\pi}\right)^n dz^1 \wedge d\ov{z^{1}} \wedge \cdots \wedge dz^n \wedge d\ov{z^{n}}.
\end{eqnarray*}
Hence
\begin{eqnarray} \nonumber
\int_X |u|_{h^m}^2 \frac{\omega^n}{n!} & \le & \frac{C}{(\log m)^2} \left( \frac{1}{a_m}\right)^n \left(1- \frac{1}{a_m}\right)^m \\ \label{eqnubound}
& \le & C (\log m)^{2n-2} m^{-\log m/2-n}.
\end{eqnarray}
%using 
%$$\left( 1 - \frac{1}{a_m} \right)^m = e^{m \log \left( 1- \frac{(\log m)^2}{m} \right)} \le e^{- (\log m)^2/2} = m^{-(\log m)/2}.$$ 
Now set $s_{m, x_0} = \eta\left( \frac{a_m |z|^2}{4} \right) v^m - u.$
Since 
$ \int_X |u|^2_{h^m} e^{-\psi} \frac{\omega^n}{n!} < \infty,$
we have $u(z) = \textrm{O}(|z|^2)$ by the definition of $\psi$. Hence $|s_{m, x_0}|^2_{h^m}(x_0) =1$.  Calculate
\begin{eqnarray} \nonumber
\int_X |s_{m, x_0}|^2_{h^m} \frac{\omega^n}{n!} & = & \int_X \left| \eta \left( \frac{a_m |z|^2}{4} \right) \right|^2 |v^m|^2_{h^m} \frac{\omega^n}{n!} + \int_X |u|^2_{h^m} \frac{\omega^n}{n!}\\ && \mbox{} + 2 \textrm{Re} \left( \int_X \left\langle \eta \left( \frac{a_m |z|^2}{4} \right) v^m, u \right\rangle_{h^m} \frac{\omega^n}{n!} \right).  \label{eqns111}
\end{eqnarray} 
From (\ref{eqnubound}) the last two terms are $\textrm{O}(m^{-q})$ for any $q$.  For the first term, observe that
\begin{equation} \label{eqns11}
\int_{|z|^2 \le 2/a_m} |v^m|^2_{h^m} \frac{\omega^n}{n!} \le \int_X \left| \eta \left( \frac{a_m |z|^2}{4} \right) \right|^2 |v^m|^2_{h^m} \frac{\omega^n}{n!} \le \int_{|z|^2 \le 4/a_m} |v^m|^2_{h^m} \frac{\omega^n}{n!}.
\end{equation}
We will make use of the following elementary lemma:

\begin{lemma} \label{lemmaintegral}
Fix $b>0$ and $q>0$.  Then
$$\left( \frac{\sqrt{-1}}{2\pi} \right)^n \int_{|z|^2 \le b/a_m} (1-|z|^2)^m dz^1 \wedge d{\ov{z^1}} \wedge \cdots \wedge dz^n \wedge d{\ov{z^n}} =  \frac{m!}{(m+n)!} + \emph{O}(m^{-q}),$$
where the term $\emph{O}(m^{-q})$ depends only on $b$, $q$ and $n$.
\end{lemma}
Then for any fixed $b>0$ and for $m$ sufficiently large,
\begin{eqnarray*}
\int_{|z|^2 \le b/a_m} |v^m|^2_{h^m} \frac{\omega^n}{n!}  & = &  \frac{m!}{(m+n)!}(1+\textrm{O}(m^{-1})).
\end{eqnarray*}
From (\ref{eqns111}) and (\ref{eqns11}) we obtain
$$\int_X |s_{m, x_0}|^2_{h^m} \frac{\omega^n}{n!} =  \frac{m!}{(m+n)!}(1+\textrm{O}(m^{-1})),$$
as required.
\end{proof1}

Assume now that  $(X, \omega)$ is a K\"ahler orbifold of complex dimension $n$ with isolated orbifold singularities at points $p_1, \ldots, p_k$.   Then for each $x \in \{ p_1, \ldots, p_k \}$ there is an open neighborhood $V_x \subset X$ containing $x$, a finite subgroup $G_x \in \textrm{U}(n)$, a $G_x$-invariant
 open neighborhood $\tilde{V}_x$ of the origin in $\mathbb{C}^n$ and a projection map $\pi_x : \tilde{V}_x \rightarrow   \tilde{V}_x/G_x \cong V_x$ with $\pi_x(0)=x$.  Let $h$ be an orbifold Hermitian metric on an orbifold line bundle $L$ with $ - \ddbar \log h = \omega$.  We show the following.

\begin{lemma} \label{lemmaorbifoldsection}
There exists $m_1>0$ depending only on $X$, $L$ and $h$ such that for each $m \ge m_1$ the following holds.   Let $x_0 \in X$.  Then
\begin{enumerate}
\item[(i)] If $x_0$
satisfies \begin{equation*} \label{eqncondition}
(d(x_0, p_i))^2 \le \frac{1}{m^2}, \quad \textrm{for some } i \in \{1, \ldots, k \},
\end{equation*}
where  $d( \ , \ )$ is the distance function on $X$ with respect to $\omega$,
then  there is a global holomorphic section $s_{m, x_0}$ of $L^m$ satisfying $|s_{m,x_0}|^2_{h^m} (x_0)  =  1$ and
\begin{eqnarray*}
%|s_{m,x_0}|^2_{h^m} (x_0) & = & 1, \\
 \frac{(m+n)!}{m!} \| s_{m, x_0} \|^2_{L^2(h^m)} & = &   \frac{1}{|G_{p_i}|} + \emph{O}(m^{-1}).
\end{eqnarray*}
\item[(ii)]  If $x_0$ satisfies
$$ \frac{1}{m^2} < (d(x_0, p_i))^2 < \frac{1}{64a_m}, \quad \textrm{for some } i \in \{1, \ldots, k \},$$
where $a_m = m/(\log m)^2$, then  there is a global holomorphic section $s_{m, x_0}$ of $L^m$ satisfying $|s_{m,x_0}|^2_{h^m} (x_0)  =  1$ and
%\begin{eqnarray*}
%|s_{m,x_0}|^2_{h^m} (x_0) & = & 1, 
%\end{eqnarray*}
\begin{eqnarray*}
 \frac{1}{C_1}\left(1  + \emph{O}(m^{-1})\right) & \le & \frac{(m+n)!}{m!} \| s_{m, x_0} \|^2_{L^2(h^m)} \le   C_2 \left(1+ \emph{O}(m^{-1}) \right),
\end{eqnarray*}
for positive constants $C_1$ and $C_2$ depending only on $|G_{p_i}|$.

\item[(iii)] If $x_0$ satisfies
$$ (d(x_0, p_i))^2 \ge \frac{1}{64a_m}, \quad \textrm{for every } i \in \{1, \ldots, k \},$$
then there exists a global holomorphic section
 $s_{m, x_0}$ of $L^m$ satisfying $|s_{m,x_0}|^2_{h^m} (x_0)  =  1$ and
\begin{eqnarray*}
%|s_{m,x_0}|^2_{h^m} (x_0) & = & 1, \\
 \frac{(m+n)!}{m!} \| s_{m, x_0} \|^2_{L^2(h^m)} & = &   1 + \emph{O}(m^{-1}).
\end{eqnarray*}
\end{enumerate}
\end{lemma}
\begin{proof1}
We will choose $m_1\gg 1$ later in the proof, depending only on $X$, $L$ and $h$.
Let $m \ge m_1$.   
For (i), 
assume first that $x_0 \in X - \{ p_1, \ldots, p_k \}$.  We may assume without loss of generality that $$(d(x_0, p_1))^2 \le \frac{1}{m^2},$$
and $d(x_0, p_i) \ge c >0$ for $i=2, \ldots, k$ for some uniform $c$.  

Dropping the subscript $p_1$, we have a uniformizing coordinate system $\pi : \tilde{V} \rightarrow \tilde{V}/G \cong V$ centered at $p_1 \in V$, so that at $0\in \tilde{V}$, the metric $g$ is the identity.  The metric is $G$-invariant and smooth in $\tilde{V}$ and has vanishing first derivatives at the origin.  Set $|G|=l$.  Since the singularity is isolated, the only fixed point of the action is $0 \in \tilde{V}$, and so $\tilde{V}- \{ 0 \} $ is a $l$-fold cover of $V- \{ p_1\}$.  The preimage of $x_0$ under the map $\pi$ consists of $l$ distinct points which we will write as $\tilde{x}_1, \ldots, \tilde{x}_l \in \tilde{V} \subset \mathbb{C}^n$.  We may assume that
$$0 < | \tilde{x}_1 |^2 = |\tilde{x}_2 |^2 = \cdots = | \tilde{x}_l|^2 < \frac{2}{m^2}.$$
Let $\eta$ and $\psi$ be the functions defined earlier in the smooth case, and  let $\tilde{\psi}$ be a weight function on $U$ given by
$$\tilde{\psi}(z) = \sum_{i=1}^l \psi(z-\tilde{x}_i).$$
Observe that $\tilde{\psi}$ is $G$-invariant, since $\psi(z)$ is a function of $|z|^2$ only.  Hence $\tilde{\psi}$ can be regarded as a smooth function on $X$ in the orbifold sense.  Note also that $\tilde{\psi}$ is non-positive everywhere.  We have
$$\sqrt{-1} \partial \overline{\partial} \tilde{\psi}(z) \ge - C(n+2) a_m \omega (z), \quad \textrm{for } z\in \tilde{V} - \{ x_1, \ldots, x_l \}.$$
Hence for sufficiently large $m$, with $\tilde{\psi}_j$ approximating $\tilde{\psi}$ as before,
$$\sqrt{-1} \partial \overline{\partial} \tilde{\psi}_j - \ddbar \log h^m + \textrm{Ric}(\omega) \ge \frac{m}{C} \omega.$$
Let $v$ be a local orbifold holomorphic section of $L$.  We may assume without loss of generality that  $|v|^2_{h}(p_1)=1$.  Pulling back to $\tilde{V}$ we have $(\partial |v|^2_{h})(0)=0$.  Let $w$ be the smooth local section of $L^m$ defined in $\tilde{V}$ by
$$w =\frac{1}{l} \sum_{i=1}^l \overline{\partial} \left( \eta \left( \frac{a_m |z-\tilde{x}_i|^2}{4} \right) \right) v^m.$$
Notice that $w$ is $G$-invariant.  If $m_1$ is sufficiently large then 
$$|\tilde{x}_i - \tilde{x}_j|^2 \le \frac{8}{m^2} \le \frac{1}{16 a_m},$$
for all $i,j$ and it follows that $w$ vanishes identically in the regions
$$\left\{ |z-\tilde{x}_i|^2 \le \frac{1}{4a_m} \right\} \quad \textrm{and} \quad  \left\{ |z-\tilde{x}_i|^2 \ge \frac{16}{a_m} \right\},$$
for all $i$.  In addition, $w$ vanishes in the regions
$$\left\{ |z|^2  \le \frac{1}{4a_m} \right\} \quad \textrm{and} \quad \left\{ |z|^2 \ge \frac{16}{a_m} \right\}.$$
It follows that $w$ descends to a smooth global orbifold section of $L^m$ and 
$\tilde{\psi}$ is uniformly bounded whenever $w$ is not identically zero.  Hence in $V$,
$$|w|^2_{h^m} e^{-\tilde{\psi}} \le C a_m |v^m|^2_{h^m}.$$
Then
$$\int_X |w|^2_{h^m} e^{-\tilde{\psi}} \frac{\omega^n}{n!} < \infty,$$
and we can apply the orbifold version of H\"ormander's estimates to obtain a smooth global orbifold section $u$ of $L^m$ satisfying $\overline{\partial} u = w$ and
$$\int_X |u|^2_{h^m} e^{-\tilde{\psi}} \frac{{\omega}^n}{n!} \le \frac{C}{m} \int_X |w|^2_{h^m} e^{-\tilde{\psi}} \frac{\omega^n}{n!} \le \frac{C}{(\log m)^2} \int_{1/4a_m \le |z|^2 \le 16/a_m} |v^m|^2_{h^m} \frac{\omega^n}{n!}.$$
We can write $|v|^2_{h}(z) = 1 - |z|^2 + \textrm{O}(|z|^3)$ and it follows that, by a similar argument as in the smooth case,
$$\int_X |u|^2_{h^m} \frac{\omega^n}{n!} \le C (\log m)^{2n-2} m^{-\log m/2 - n}.$$
Now set
$$s_{m, x_0}(z)= \frac{1}{l} \sum_{i=1}^l \eta \left( \frac{a_m |z-\tilde{x}_i|^2}{4} \right) v^m(z) - u(z),$$
so that $s_{m,x_0}$ is a global holomorphic orbifold section of $L^m$.  Notice that since
$$\int_X |u|^2_{h^m} e^{-\tilde{\psi}} \frac{\omega^n}{n!} < \infty,$$
it follows from the definition of $\tilde{\psi}$ that $u(z-\tilde{x}_i) = \textrm{O}(|z-\tilde{x}_i|^2)$ for each $i$.  Hence $u(x_0)=0$ and since $|s_{m,x_0}|_{h^m}^2(x_0) =  |v^m|^2_{h^m}(x_0)$, we have
\begin{equation}\label{lemmas2}
1   \ge  |s_{m,x_0}|_{h^m}^2(x_0)   \ge (1- \frac{2}{m^2})^m 
 =  e^{m \log (1 - 2/m^2)} 
 =  1 - \textrm{O}(m^{-q}),
\end{equation}
for any $q$.
Calculate, remembering that $\tilde{V}-\{ 0 \}$ is an $l$-fold cover of $V - \{ p_1 \}$,
\begin{eqnarray} \nonumber
\int_X |s_{m,x_0}|^2_{h^m} \frac{\omega^n}{n!} & = & \frac{1}{l} \int_{\tilde{V}} \left| \frac{1}{l} \sum_{i=1}^l \left( \eta \left( \frac{a_m |z-\tilde{x}_i|^2}{4} \right) \right) \right|^2 |v^m|^2_{h^m} \frac{\omega^n}{n!} \\ \nonumber && \mbox{} +  \frac{2}{l}\textrm{Re} \left( \int_{\tilde{V}} \left\langle \frac{1}{l} \sum_{i=1}^l  \eta \left( \frac{a_m |z-\tilde{x}_i|^2}{4} \right) v^m , u \right\rangle_{h^m} \frac{\omega^n}{n!} \right)  \\ && \mbox{}
+ \frac{1}{l} \int_X |u|^2_{h^m} \frac{\omega^n}{n!}.\qquad \label{eqns}
\end{eqnarray} 
The last two terms are $\textrm{O}(m^{-q})$ for any $q$.  For the first term, observe that
\begin{eqnarray}  \nonumber
\frac{1}{l}\int_{|z|^2 \le 1/4a_m} |v^m|^2_{h^m} \frac{\omega^n}{n!} & \le & \frac{1}{l} \int_{\tilde{V}} \left| \frac{1}{l} \sum_{i=1}^l \left( \eta \left( \frac{a_m |z-\tilde{x}_i|^2}{4} \right) \right) \right|^2 |v^m|^2_{h^m} \frac{\omega^n}{n!}  \\ & \le & \frac{1}{l} \int_{|z|^2 \le 16/a_m} |v^m|^2_{h^m} \frac{\omega^n}{n!}. \label{eqns1}
\end{eqnarray}
From Lemma \ref{lemmaintegral}  we obtain
$$\int_X |s_{m,x_0}|_{h^m}^2 \frac{\omega^n}{n!} = \frac{m!}{l(m+n)!} (1+\textrm{O}(m^{-1})).$$
Then from (\ref{lemmas2}) we obtain the required section by rescaling.  

The case when $x_0$ is one of the singular points $p_i$ is easier, since we can take the weight function to be $\psi$, which is of course $G$-invariant.  The proof follows as in the smooth case, except that a factor of $l$ arises when estimating the integral of $|s_{m, x_0}|^2_{h^m}$.

We now consider case (ii).   We divide this into two parts:
\begin{enumerate}
\item[(a)] $\displaystyle{\frac{1}{m^2} < (d(x_0, p_1))^2 < \frac{A}{m}}$;
\item[(b)] $\displaystyle{\frac{A}{m} \le (d(x_0, p_1))^2 < \frac{1}{64a_m}}$;
\end{enumerate}
for a constant $A \gg 0$ depending on $l$ to be determined later.  

For (a), we can use almost the same argument as in (i).  The only difference is that (\ref{lemmas2}) becomes
$$ 1 \ge |s_{m, x_0}|^2_{h^m}(x_0) \ge c>0,$$
for a constant $c$ depending on $A$.  The required estimate follows after scaling $s_{m, x_0}$.

For (b) we argue as follows.
Using the  notation above, we work in the coordinate patch $\tilde{V}$ and
consider the same weight function $\tilde{\psi}$.  Let $v_1$ be a local holomorphic section of $\pi^*L$ over $\tilde{V}$  with the property that
$$|v_1|^2_h(z) = 1 - |z-\tilde{x}_1|^2 + \textrm{O}(|z-\tilde{x}_1|^3).$$
Observe that $v_1$ is not $G$-invariant.  Writing the elements of $G$ as $\gamma_1, \ldots, \gamma_l$ with $\gamma_1$ the identity element, we set $v_i = \gamma_i^* v_1$, so that $$|v_i|^2_h(z) = 1 - |z-\tilde{x}_i|^2 + \textrm{O}(|z-\tilde{x}_i|^3).$$
Now define a local $G$-invariant section $\hat{w}$ of $L^m$  over $\tilde{V}$ by
$$\hat{w}(z) =\frac{1}{l} \sum_{i=1}^l \overline{\partial} \left( \eta \left( \frac{a_m |z-\tilde{x}_i|^2}{4} \right) \right) v_i^m(z).$$  We have
\begin{equation} \label{eqnxixj}
\frac{A}{4m} \le | \tilde{x}_i - \tilde{x}_j|^2 \le \frac{1}{16a_m},
\end{equation}
and by a similar argument as in case (i),
$$\int_X |\hat{w}|^2_{h^m} e^{-\tilde{\psi}} \frac{\omega^n}{n!} < \infty.$$
Hence we can obtain a smooth global orbifold section $\hat{u}$ of $L^m$ satisfying $\ov{\partial} \hat{u} = \hat{w}$ and
$$\int_X |\hat{u}|^2 \frac{\omega^n}{n!} \le C ( \log m)^{2n-2} m^{-\log m/2 -n}.$$
Let $s_{m, x_0}$ be the global holomorphic orbifold section of $L^m$ given by
$$s_{m, x_0}(z) =\frac{1}{l} \sum_{i=1}^l   \eta \left( \frac{a_m |z-\tilde{x}_i|^2}{4} \right)  v_i^m(z) - \hat{u}(z).$$
Notice that $\hat{u}(x_0)=0$.   For $i \neq j$ we have
$|v_i^m|_{h^m}^2(\tilde{x}_j) \le e^{-A/8}$, using (\ref{eqnxixj}),
and if $A$ is sufficiently large depending only on $l$ it follows that 
\begin{equation} \label{eqniia}
1 \ge  |s_{m, x_0}|^2_{h^m}(x_0) \ge c>0,
\end{equation}
for $c$ depending only on $l$.  Now
\begin{eqnarray} \nonumber
\int_X |s_{m,x_0}|^2_{h^m} \frac{\omega^n}{n!} & = & \frac{1}{l} \int_{\tilde{V}} \left| \frac{1}{l} \sum_{i=1}^l  \eta \left( \frac{a_m |z-\tilde{x}_i|^2}{4}  \right) v_i^m \right|^2_{h^m} \frac{\omega^n}{n!} \\ \nonumber && \mbox{} +  \frac{2}{l}\textrm{Re} \left( \int_{\tilde{V}} \left\langle \frac{1}{l} \sum_{i=1}^l  \eta \left( \frac{a_m |z-\tilde{x}_i|^2}{4} \right) v_i^m, \hat{u} \right\rangle_{h^m} \frac{\omega^n}{n!} \right)  \\ && \mbox{}
+ \frac{1}{l} \int_X |\hat{u}|^2_{h^m} \frac{\omega^n}{n!}.\qquad \label{eqnsiii}
\end{eqnarray} 
As before, the last two terms are $\textrm{O}(m^{-q})$ for any $q$, and 
\begin{eqnarray}  \nonumber
 c'\frac{m!}{(m+n!)}(1+ \textrm{O}(m^{-1})) & \le & \frac{1}{l} \int_{\tilde{V}} \left| \frac{1}{l} \sum_{i=1}^l  \eta \left( \frac{a_m |z-\tilde{x}_i|^2}{4} \right) v_i^m \right|^2_{h^m} \frac{\omega^n}{n!}  \\ & \le &\frac{m! }{l(m+n!)}(1+ \textrm{O}(m^{-1})), \label{eqns1iii}
\end{eqnarray}
for a constant $c'>0$ depending only on $l$.  Combining (\ref{eqniia}), (\ref{eqnsiii}) and (\ref{eqns1iii}) completes part (b) of (ii).  

For case (iii) we can avoid the singularities using the same argument as in the smooth case.
\end{proof1}

\begin{remark}
The result of Lemma \ref{lemmaorbifoldsection}.(ii) is clearly not sharp.  It would be interesting to know what the optimal estimates are in this case. \end{remark}

\setcounter{equation}{0}
\setcounter{lemma}{0}
\addtocounter{section}{1}
\bigskip
\noindent
{\bf 3.  Convergence of Tsuji's iteration}
\bigskip

In this section we give a proof of Theorem 1.  We begin with a simple and well-known observation, which we will be useful later.  Let $\mathcal{X}$ be any set, and let $H$ be a finite dimensional vector subspace of the vector space of functions from $\mathcal{X}$ to $\mathbb{C}$.  Suppose that $H$ is equipped with an inner product $\langle \ , \ \rangle_H$.  For any orthonormal basis $(v_0, \ldots, v_N)$ of $H$, define a function $\rho: \mathcal{X} \rightarrow \mathbb{R}$ by
$\rho (x) = \sum_{i=0}^N |v_i|^2(x).$
Note that the function $\rho$ is independent of the choice of orthonormal basis.  Now fix $x \in \mathcal{X}$.  Then it is possible to choose an orthonormal basis $(v_0, \ldots, v_N)$ such that 
$v_i (x) = 0$ for $i=1, 2, \ldots, N$.
The observation is that
\begin{equation} \label{eqnobservation}
\rho(x) = \sup \left\{ |v|^2(x) \ \big| \ v \in H, \| v \|_H =1  \right\} = | v_0|^2(x).
\end{equation}

We now turn to the proof of Theorem 1.   Notice that, in addition to a Hermitian metric $h_{m}$ on $K_X^m$ for each $m\ge m_0$, we have defined by (\ref{eqnT}) an inner product $\langle \ , \ \rangle_{T_{m}}$ on $\HZ(X, K_X^{m})$ for each $m \ge m_0 +1$.  Also, from the Hermitian metric $\hke$ on $K_X$ we have the $L^2$ inner product on $\HZ (X, K_X^m)$ given by $\hke^m$ and $\frac{1}{n!}\omega_{\textrm{KE}}^n$.  We will denote this inner product simply by $\langle \cdot, \cdot \rangle_{\textrm{KE}}$.

We will use Lemma \ref{lemmapeaksections} on the existence of peak sections to prove the following:

\begin{lemma} \label{lemmaineq}
Let $m_1$ be the constant of Lemma \ref{lemmapeaksections} for $L=K_X$, $h=h_{\emph{KE}}$. Assume $m_1 \ge m_0$. There exists $A$ depending only on $X$ such that for all $m \ge m_1$,
\begin{equation} \label{eqnhmineq}
\inf_X \left( \frac{h_{m_1}}{h^{m_1}_{\emph{KE}}} \right)   h_{\emph{KE}}^m \prod_{k=m_1}^m \left( 1- \frac{A}{k} \right) \le h_m \le \sup_X \left( \frac{h_{m_1}}{h^{m_1}_{\emph{KE}}} \right) h_{\emph{KE}}^{m} \prod_{k=m_1}^m \left( 1 + \frac{A}{k} \right).
\end{equation}
\end{lemma}
\begin{proof1}
In the course of this proof, the constant $A$ may change from line to line.
We will prove the upper bound on $h_m$ first.  We use induction.  Obviously, the inequality holds for $m=m_1$.  Let 
$$C_m = \sup_X \left( \frac{h_{m_1}}{h^{m_1}_{\emph{KE}}} \right) \prod_{k=m_1}^m \left( 1 + \frac{A}{k} \right),$$
and assume that $h _m \le C_m  \hke^{m}.$  Notice that for any section $t$ of $\HZ (X, K_X^{m+1})$,
$$ \left\| t \right\|^2_{T_{m+1}} \le C_m \left\| t \right\|_{\textrm{KE}}^2.$$
Fix a point $x_0 \in X$.  Let $s_{m+1,x_0} \in \HZ (X, K_X^{m+1})$ be a peak section as constructed in Lemma \ref{lemmapeaksections}.  
Then calculate, from the definition of $h_{m+1}$,
\begin{eqnarray*}
\frac{(m+1)!}{(m+n+1)!} h_{m+1}( x_0) & \le & \frac{ \left\| s_{m+1, x_0} \right\|^2_{T_{m+1}}}{  \left|s_{m+1, x_0}\right|_{\hke^{m+1}}^2(x_0)}\hke^{m+1}(x_0) \\
& \le & C_m \left\| s_{m+1, x_0} \right\|_{KE}^2 \, \hke^{m+1}(x_0) \\
& = & C_m \frac{(m+1)!}{(m+n+1)!} \left(1+ \textrm{O}\left( \frac{1}{m+1} \right) \right) \hke^{m+1}(x_0), \\
\end{eqnarray*}
and it follows that
$$h_{m+1} \le \sup_X \left( \frac{h_{m_1}}{h^{m_1}_{\emph{KE}}}  \right) h_{\emph{KE}}^{m+1}\prod_{k=m_1}^{m+1} \left( 1 + \frac{A}{k} \right).$$
We turn now to the lower bound for $h_m$.  Again we use induction and assume that $h_m \ge D_m \hke^m$ for 
$$D_m = \inf_X \left( \frac{h_{m_1}}{h^{m_1}_{\emph{KE}}} \right) \prod_{k=m_1}^m \left( 1 - \frac{A}{k} \right).$$
Fix $x_0$ in $X$.
Let $(\sigma_{m+1}^{(0)}, \ldots, \sigma_{m+1}^{(N_{m+1})})$ be an orthonormal basis of $\HZ(X, K_X^{m+1})$ with respect to the inner product $\langle \ , \ \rangle_{T_{m+1}}$.  We may assume that
$$\sigma_{m+1}^{(i)}(x_0)=0 \quad \textrm{for } i=1, \ldots, N_{m+1}.$$
Observe that $1 = \| \sigma_{m+1}^{(0)} \|^2_{T_{m+1}} \ge D_m \| \sigma_{m+1}^{(0)} \|_{\textrm{KE}}^2$.  Then
\begin{eqnarray} \nonumber
\frac{(m+1)!}{(m+n+1)!} h_{m+1}( x_0) & = & \hke^{m+1}(x_0)\frac{1}{ | \sigma^{(0)}_{m+1} |^2_{\hke^{m+1}}(x_0)} \\
& \ge & D_m \hke^{m+1}(x_0) \frac{ \| \sigma^{(0)}_{m+1} \|_{\textrm{KE}}^2}{| \sigma^{(0)}_{m+1} |^2_{\hke^{m+1}}(x_0)}.
\end{eqnarray}
Now let $( \tau_{m+1}^{(0)}, \ldots, \tau_{m+1}^{(N_{m+1})} )$ be an orthonormal basis of $\HZ(X, K_X^{m+1})$ with respect to the inner product $\langle \ , \ \rangle_{\textrm{KE}}$.  As before we may assume that
$$\tau_{m+1}^{(i)}(x_0)=0 \quad \textrm{for } i=1, \ldots, N_{m+1}.$$
Then it follows that if $t$ is any section of $\HZ(X, K_X^{m+1})$, we have
\begin{eqnarray}
\frac{ |t|^2_{\hke^{m+1}} (x_0)}{ \| t \|^2_{\textrm{KE}}} \le { | \tau^{(0)}_{m+1} |^2_{\hke^{m+1}} (x_0)}. 
\end{eqnarray}
Hence
\begin{eqnarray} 
\frac{(m+1)!}{(m+n+1)!} h_{m+1}( x_0) & \ge  & D_m \hke^{m+1}(x_0) \frac{ 1}{| \tau^{(0)}_{m+1} |^2_{\hke^{m+1}}(x_0)}. \label{eqnhlowerbound1}
\end{eqnarray}
We consider again the peak section $s_{ m+1, x_0}$.  Define real numbers $a_0, \ldots, a_{N_{m+1}}$ by
$$ s_{m+1, x_0} = \sum_{i=0}^{N_{m+1}} a_i \tau_{m+1}^{(i)}.$$
Then, using the second part of Lemma \ref{lemmapeaksections},
\begin{eqnarray*}
\sum_{i=1}^{N_{m+1}} a_i^2 & = &\left\langle s_{m+1, x_0} , \sum_{i=1}^{N_{m+1}} a_i \tau_{m+1}^{(i)} \right\rangle_{\textrm{KE}} \\
& \le & \frac{A}{m+1} \| s_{m+1, x_0} \|_{\textrm{KE}} \left( \sum_{i=1}^{N_{m+1}} a_i^2
\right)^{1/2},
\end{eqnarray*}
and so 
$$ \sum_{i=1}^{N_{m+1}} a_i^2 \le A \frac{ \| s_{m+1, x_0} \|^2_{\textrm{KE}}}{(m+1)^2}.$$
Now notice that
\begin{eqnarray*}
a_0^2 & = & \| s_{m+1, x_0} \|^2_{\textrm{KE}} - \sum_{i=1}^{N_{m+1}} a_i^2 \\
& \ge & \| s_{m+1, x_0} \|^2_{\textrm{KE}} \left(1 - A \frac{1}{(m+1)^2} \right)
\end{eqnarray*}
Now since $|s_{m+1, x_0}|^2_{\hke^{m+1}}(x_0)=1$, we have $|\tau^{(0)}_{m+1}|^2_{\hke^{m+1}}(x_0)=1/a_0^2$.  Then from (\ref{eqnhlowerbound1}) we have
\begin{eqnarray*} 
h_{m+1}( x_0) & \ge  & D_m \hke^{m+1}(x_0)  \left(1 -  \textrm{O}\left( \frac{1}{m+1}\right) \right),
\end{eqnarray*}
and the required lower bound follows.
\end{proof1}

From this, we can prove Theorem 1.

\bigskip
\noindent
{\bf Proof of Theorem 1} \ Raise (\ref{eqnhmineq}) to the power $1/m$. 
For the upper bound of Theorem 1, observe that
\begin{eqnarray*}
\log \left( \prod_{m_1}^m \left(1+\frac{A}{k} \right) \right)^{1/m} & \le & \frac{1}{m} \sum_{k=1}^m \log \left(1 + \frac{A}{k} \right) \\
& \le & \frac{C}{m} \sum_{k=1}^m \frac{1}{k} \\
& \le & C \frac{\log m}{m},
\end{eqnarray*}
for a constant $C$ depending only on $A$.
The lower bound follows similarly. \hfill$\square$\medskip

\setcounter{equation}{0}
\setcounter{lemma}{0}
\addtocounter{section}{1}
\bigskip
\noindent
{\bf 4.  The modified iteration}
\bigskip

We give a proof of Theorem 2.  We omit the proof of Theorem 3, since it is simpler and follows along the same lines.  
We first prove a convergence result on the minimal surface $\Xmin$ of general type.   Consider a Hermitian metric $h_{m_0}$ on $K_{\Xmin}^{m_0}$ and write $\beta = |\SM|^2_{\hC}$, for $|\SM|^2_{\hC}$ as in the introduction.  Recall that $\hmin$ is the Hermitian metric on $K_{\Xmin}$ given by $f^*\hke$, for $\hke$ the Hermitian metric on $K_{\Xcan}$ corresponding to the K\"ahler-Einstein metric $\omega_{\textrm{KE}}$. Consider the sequence of metrics $h_{m, \ep} = h_{m, \ep}(\beta, h_{m_0})$ on $K_{\Xmin}^m$.  

\begin{theorem} \label{theoremmin}
For every sequence $\ep_j \rightarrow 0$, 
$$\limsup_{m \rightarrow \infty} h_{m, \ep_j}^{1/m} \rightarrow h_{\emph{min}}, \ \ \textrm{as } j \rightarrow \infty,$$
almost everywhere on $X_{\emph{min}}$.
\end{theorem}

To prove this, we will need two lemmas.

\begin{lemma} \label{lemmaorbifoldupperbound}
There exist $m_1>0$  and $A$ depending only on $X_{\emph{min}}$, $\beta$ and $\ep$ such that for all $m \ge  m_1$,
\begin{equation} \label{eqngentypeupperbound}
\beta^{\ep} h_{m, \ep} \le \sup_{X_{\emph{min}}} \left( \frac{h_{m_1, \ep}}{h^{m_1}_{\emph{min}}} \right) h_{\emph{min}}^{m} \prod_{k=m_1}^m \left( 1 + \frac{A}{k} \right).
\end{equation}
\end{lemma}
\begin{proof1}
We will use induction, for $m_1$ to be determined later.  Assume the inequality holds for $h_m$.  Let
$$C_m = \sup_{\Xmin} \left( \frac{h_{m_1, \ep}}{h^{m_1}_{\textrm{min}}} \right) \prod_{k=m_1}^m \left( 1 + \frac{A}{k} \right).$$
Then $\beta^{\ep} h_{m, \ep} \le C_m \hmin^m.$  It follows that $\| t \|^2_{T_{m+1, \ep}} \le C_m \| t \|_{L^2(\hmin^{m+1})}^2$ for $t$ any global section of $K_{\Xmin}^{m+1}$.  Since the inequality for $h_{m+1, \ep}$ obviously holds at points on $\mathcal{C}$, it is sufficient to prove it  at  a fixed point $y_0 \in \Xmin - \mathcal{C}$.  Write $x_0 = f(y_0) \in \Xcan$.  By Lemma \ref{lemmaorbifoldsection} there is a global holomorphic section $s_{m+1, x_0}$ of $K_{\Xcan}^{m+1}$ satisfying $|s_{m+1, x_0}|^2_{h_{\textrm{KE}}^{m+1}}(x_0) =1$ and
$$ \beta^{\ep} (y_0) \int_{\Xcan} |s_{m+1, x_0}|^2_{h_{\textrm{KE}}^{m+1}} \frac{\omega_{\textrm{KE}}^n}{n!} \le \frac{(m+1)!}{(m+n+1)!} \left( 1 + \textrm{O}(m^{-1}) \right),$$
as long as $m_1$ is chosen to be sufficiently large.  Then
\begin{eqnarray*}
\lefteqn{\frac{(m+1)!}{(m+n+1)!}  \beta^{\ep} (y_0) h_{m+1, \ep}( y_0) } \\& \le & \beta^{\ep}(y_0) \frac{ \left\| f^*(s_{m+1, x_0}) \right\|^2_{T_{m+1, \ep}}}{  \left| f ^*(s_{m+1, x_0}) \right|_{\hmin^{m+1}}^2(y_0)}\hmin^{m+1}(y_0) \\
& \le & C_m  \beta^{\ep}(y_0)   \left( \int_{\Xcan} |s_{m+1, x_0}|^2_{h_{\textrm{KE}}^{m+1}} \frac{\omega_{\textrm{KE}}^n}{n!} \right)  \hmin^{m+1}(y_0) \\
& \le & C_m \frac{(m+1)!}{(m+n+1)!} \left(1+ \textrm{O}\left( m^{-1} \right) \right) \hmin^{m+1}(y_0),
\end{eqnarray*}
and the lemma follows.
\end{proof1}

For the lower bound of $h_{m, \ep}$, we use a modification of a lemma of Tsuji \cite{Ts}:

\begin{lemma} \label{lemmaintegralbound}  
There exist constants $m_2$ and $B$ depending only on $X_{\emph{min}}$ such that for all $m\ge m_2$ and $0< \ep \le 1$,
\begin{equation} \label{eqnintegralbound}
\int_{X_{\emph{min}}} \beta^{\ep} h_{m, \ep}^{-1/m} \le V^{\frac{m-m_2+1}{m}} \left( \prod_{k=m_2}^m \left(1 + \frac{B}{k} \right) \right)^{1/m} \left( \int_{X_{\emph{min}}} \beta^{\ep} h_{m_2, \ep}^{-1/m_2} \right)^{m_2/m},
\end{equation}
for $\displaystyle{V=  \int_{X_{\emph{min}}} \frac{\omega_{\min}^n}{n!}}$.
\end{lemma}
\begin{proof1}
Set $$L_{m, \ep} = \int_{\Xmin} \beta^{\ep} h_{m,\ep}^{-1/m}$$ and $$c_m = \frac{m!}{(m+n)!} (N_m+1),$$
for $N_m+1 = \textrm{dim} \, \textrm{H}^0(\Xmin, K_{\Xmin}^m)$.  Then we claim that
\begin{equation} \label{eqnrecursive}
L_{m, \ep} \le c_m^{1/m} L_{m-1, \ep}^{(m-1)/m}.
\end{equation}
Given (\ref{eqnrecursive}), we can finish the proof of the lemma as follows.  First, by Riemann-Roch, there exist constants $m_2$ and $B$ such that for $m \ge m_2$,
$$ c_m \le V \left(1+\frac{B}{m} \right).$$
From (\ref{eqnrecursive}), arguing by induction, we have
$$L_{m, \ep} \le \left( c_m c_{m-1} \cdots c_{m_2} \right)^{1/m} L_{m_2, \ep}^{m_2/m},$$
and the inequality (\ref{eqnintegralbound}) follows immediately.
It remains to show (\ref{eqnrecursive}).  Using H\"older's inequality,
\begin{eqnarray*}
L_{m, \ep} 
& = & \int_{\Xmin} \beta^{\ep} h_{m, \ep}^{-1/m} h_{m-1, \ep}^{1/(m-1)} h_{m-1, \ep}^{-1/(m-1)} \\
& \le & \left( \int_{\Xmin} \beta^{\ep}(h_{m, \ep}^{-1/m} h_{m-1, \ep}^{1/(m-1)} )^m h^{-1/(m-1)}_{m-1, \ep} \right)^{\frac{1}{m}} \left( \int_{\Xmin} \beta^{\ep} h_{m-1, \ep}^{-1/(m-1)} \right)^{\frac{m-1}{m}} \\
& = & \left( \int_{\Xmin} \beta^{\ep} h^{-1}_{m, \ep} h_{m-1, \ep} \right)^{\frac{1}{m}} L_{m-1, \ep}^{(m-1)/m} \\
& = & \left( \frac{m!}{(m+n)!} \int_{\Xmin} \beta^{\ep} \sum_{i=0}^{N_m} \sigma_{m, \ep}^{(i)} \otimes \overline{\sigma_{m, \ep}^{(i)}} \otimes h_{m-1, \ep} \right)^{\frac{1}{m}} L_{m-1, \ep}^{(m-1)/m} \\
& = & c_m^{1/m} L_{m-1, \ep}^{(m-1)/m},
\end{eqnarray*}
and this completes the proof of the lemma.
\end{proof1}

We can now use these lemmas to prove a convergence result for the metrics $h_{m, \ep}$.

\bigskip
\noindent
{\bf Proof of Theorem \ref{theoremmin}} \ 
From Lemma \ref{lemmaorbifoldupperbound} we have 
$$h_{m, \ep}^{-1/m} \ge E(m, \ep) \hmin^{-1} \quad \textrm{on } \Xmin- \mathcal{C},$$
where $E(m, \ep) \rightarrow 1$ as $m \rightarrow \infty$.  From Lemma \ref{lemmaintegralbound}, we have
$$\int_{\Xmin} \beta^{\ep} h_{m, \ep}^{-1/m} \le V \, F(m, \ep),$$
where $F(m, \ep) \rightarrow 1$ as $m \rightarrow \infty$.   Writing $h_{\ep} = \limsup_{m \rightarrow \infty} h_{m, \ep}^{1/m}$, we have
\begin{eqnarray*}
\int_{\Xmin} \beta^{\ep} \left| \frac{h_{\ep}^{-1} - \hmin^{-1}}{\hmin^{-1}} \right| \hmin^{-1}
& = & \int_{\Xmin} \liminf_{m \rightarrow \infty} \beta^{\ep} (h_{m, \ep}^{-1/m} - E(m, \ep) \hmin^{-1}) \\
& \le & \liminf_{m \rightarrow \infty} \int_{\Xmin} \beta^{\ep} h_{m, \ep}^{-1/m} - \int_{\Xmin} \beta^{\ep}  \hmin^{-1} \\
& \le & V - \int_{\Xmin} \beta^{\ep} \hmin^{-1} \rightarrow 0,
\end{eqnarray*}
as $\ep \rightarrow 0$.  Theorem \ref{theoremmin} follows.
\hfill$\square$\medskip

Finally, we complete the proof of Theorem 2.

\bigskip
\noindent
{\bf Proof of Theorem 2} \ Using the notation given in the introduction, there is an isomorphism $\Theta : \HZ(\Xmin, K_{\Xmin}^m) \rightarrow \HZ(X, K_X^m)$ given by $\Theta (s) = \tau^* s \otimes \SE^m$.  Then, given an inner product $T_{m, \ep}$ on $\HZ(X, K_X^m)$, we can define an inner product $\hat{T}_{m, \ep}$ on $\HZ(\Xmin, K_{\Xmin}^m)$ by $\langle s, t \rangle_{\hat{T}_{m, \ep}} = \langle \Theta(s), \Theta(t) \rangle_{T_{m, \ep}}$. 

Then given an initial Hermitian metric $h_{m_0}$ on $K_X^{m_0}$ we can obtain an inner product $\langle \cdot, \cdot \rangle_{\hat{T}_{m_0+1, \ep}}$ on $K_{\Xmin}^{m_0+1}$ and hence a Hermitian metric $\hat{h}_{m_0+1, \ep}$ on $K_{\Xmin}^{m_0+1}$.  Applying the modified Tsuji iteration as in the case of Theorem \ref{theoremmin}, 
we obtain a sequence of Hermitian metrics $\hat{h}_{m, \ep}$ for $m \ge m_0+1$ on $K_{\Xmin}^m$.  From the definition of $\SE$, one can check that $h_{m, \ep} = \tau^* \hat{h}_{m, \ep} \otimes |\SE|^{-2m}$.  Indeed, assuming inductively that $h_{m, \ep} = \tau^* \hat{h}_{m, \ep} \otimes |\SE|^{-2m}$, denote by  $\langle \cdot, \cdot, \rangle'_{\hat{T}_{m+1, \ep}}$ the inner product induced by $\hat{h}_{m, \ep}$ on $\HZ(\Xmin, K_{\Xmin}^{m+1})$.  We need to show that $\langle s, t \rangle'_{\hat{T}_{m+1, \ep}} =\langle s, t \rangle_{\hat{T}_{m+1, \ep}}$ for $s, t \in \HZ(\Xmin, K_{\Xmin}^{m+1})$.  But
\begin{eqnarray*}
\langle s, t \rangle'_{\hat{T}_{m+1, \ep}} & = & \int_{\Xmin} \hat{h}_{m, \ep} \otimes s \otimes \ov{t} \\
& = & \int_X h_{m, \ep} \otimes |S_{-1}|^{2m} \otimes (\tau^* s) \otimes (\ov{\tau^* t}) \otimes |S_{-1}|^2 \\
& = & \langle \Theta(s), \Theta(t) \rangle_{T_{m+1, \ep}} \\
& = & \langle s, t \rangle_{\hat{T}_{m+1, \ep}}.
\end{eqnarray*}
Now by Theorem \ref{theoremmin}, we see that for any sequence $\ep_j \rightarrow 0$, we have $$\limsup_{m \rightarrow \infty} \hat{h}_{m, \ep_j}^{1/m} \rightarrow \hmin,$$ almost everywhere.  Theorem 2 follows immediately.  \hfill$\square$

\end{document}